\documentclass[12pt]{article}
\usepackage{amsmath,amssymb,amsbsy,amsfonts,amsthm,latexsym,
                     amsopn,amstext,amsxtra,euscript,amscd}

\def\bbbc{{\mathchoice {\setbox0=\hbox{$\displaystyle\rm C$}\hbox{\hbox
to0pt{\kern0.4\wd0\vrule height0.9\ht0\hss}\box0}}
{\setbox0=\hbox{$\textstyle\rm C$}\hbox{\hbox
to0pt{\kern0.4\wd0\vrule height0.9\ht0\hss}\box0}}
{\setbox0=\hbox{$\scriptstyle\rm C$}\hbox{\hbox
to0pt{\kern0.4\wd0\vrule height0.9\ht0\hss}\box0}}
{\setbox0=\hbox{$\scriptscriptstyle\rm C$}\hbox{\hbox
to0pt{\kern0.4\wd0\vrule height0.9\ht0\hss}\box0}}}}
\def\bbbq{{\mathchoice {\setbox0=\hbox{$\displaystyle\rm
Q$}\hbox{\raise
0.15\ht0\hbox to0pt{\kern0.4\wd0\vrule height0.8\ht0\hss}\box0}}
{\setbox0=\hbox{$\textstyle\rm Q$}\hbox{\raise
0.15\ht0\hbox to0pt{\kern0.4\wd0\vrule height0.8\ht0\hss}\box0}}
{\setbox0=\hbox{$\scriptstyle\rm Q$}\hbox{\raise
0.15\ht0\hbox to0pt{\kern0.4\wd0\vrule height0.7\ht0\hss}\box0}}
{\setbox0=\hbox{$\scriptscriptstyle\rm Q$}\hbox{\raise
0.15\ht0\hbox to0pt{\kern0.4\wd0\vrule height0.7\ht0\hss}\box0}}}}
\def\bbbt{{\mathchoice {\setbox0=\hbox{$\displaystyle\rm
T$}\hbox{\hbox to0pt{\kern0.3\wd0\vrule height0.9\ht0\hss}\box0}}
{\setbox0=\hbox{$\textstyle\rm T$}\hbox{\hbox
to0pt{\kern0.3\wd0\vrule height0.9\ht0\hss}\box0}}
{\setbox0=\hbox{$\scriptstyle\rm T$}\hbox{\hbox
to0pt{\kern0.3\wd0\vrule height0.9\ht0\hss}\box0}}
{\setbox0=\hbox{$\scriptscriptstyle\rm T$}\hbox{\hbox
to0pt{\kern0.3\wd0\vrule height0.9\ht0\hss}\box0}}}}
\def\bbbs{{\mathchoice
{\setbox0=\hbox{$\displaystyle     \rm S$}\hbox{\raise0.5\ht0\hbox
to0pt{\kern0.35\wd0\vrule height0.45\ht0\hss}\hbox
to0pt{\kern0.55\wd0\vrule height0.5\ht0\hss}\box0}}
{\setbox0=\hbox{$\textstyle        \rm S$}\hbox{\raise0.5\ht0\hbox
to0pt{\kern0.35\wd0\vrule height0.45\ht0\hss}\hbox
to0pt{\kern0.55\wd0\vrule height0.5\ht0\hss}\box0}}
{\setbox0=\hbox{$\scriptstyle      \rm S$}\hbox{\raise0.5\ht0\hbox
to0pt{\kern0.35\wd0\vrule height0.45\ht0\hss}\raise0.05\ht0\hbox
to0pt{\kern0.5\wd0\vrule height0.45\ht0\hss}\box0}}
{\setbox0=\hbox{$\scriptscriptstyle\rm S$}\hbox{\raise0.5\ht0\hbox
to0pt{\kern0.4\wd0\vrule height0.45\ht0\hss}\raise0.05\ht0\hbox
to0pt{\kern0.55\wd0\vrule height0.45\ht0\hss}\box0}}}}
\def\bbbz{{\mathchoice {\hbox{$\sf\textstyle Z\kern-0.4em Z$}}
{\hbox{$\sf\textstyle Z\kern-0.4em Z$}}
{\hbox{$\sf\scriptstyle Z\kern-0.3em Z$}}
{\hbox{$\sf\scriptscriptstyle Z\kern-0.2em Z$}}}}

\def\cU{\mathcal U}

\def\e{{\mathbf{\,e}}}
\def\el{{\mathbf{\,e}}_\ell}
\def\eq{{\mathbf{\,e}}_q}
\def\eqd{{\mathbf{\,e}}_{q_d}}

\def\Z{\mathbb{Z}}

\renewcommand{\vec}[1]{\mathbf{#1}}

\def\vec#1{\mathbf{#1}}
\def\inv#1{\overline{#1}}

\def\bbbc{{\mathrm I\!C}}

\def\bbbq{{\mathbb Q}}

\def \Z{{\bbbz}}

\newtheorem{thm}{Theorem}% [section]
\newtheorem{lem}[thm]{Lemma}

\begin{document}

\def\\{\cr}
\def\({\left(}
\def\){\right)}
\def\fl#1{\left\lfloor#1\right\rfloor}
\def\rf#1{\left\lceil#1\right\rceil}

\title{On a Generalised Lehmer Problem for Arbitrary 
Powers}

\author{ 
{\sc Igor E.~Shparlinski} \\
{Department of Computing}\\
{Macquarie University} \\
{Sydney, NSW 2109, Australia} \\
{\tt igor@ics.mq.edu.au}}

\maketitle

\begin{abstract}
We consider a generalisation of the classical Lehmer problem
about the parity distribution of an integer and its modular 
inverse. We use some known estimates of exponential
sums to study a more general question of 
simultaneous distribution 
of the residues of any fixed number of negative and positive 
powers of integers in prescribed arithmetic 
progressions. In particular,  we improve and generalise  a recent
result of Y.~Yi and W.~Zhang. 
\end{abstract}

%%%%%%%%%%%Introducere%%%%%%%%%%%

\section{Introduction} 

Given modulus $q\ge 2$, we denote by $\cU_q$ 
the set 
$$
\cU_q = \{n \ : \ 1 \le n < q, \ \gcd(n,q)=1\}.
$$
In particular, $\# \cU_q = \varphi(q)$, the Euler function.

For $n \in \cU_q$ we use $\inv{n}$ to denote 
the modular inverse
of 
$n$, that is, the unique integer $\inv{n}\in \cU_q$ 
with  $n \inv{n} \equiv 1 \pmod q$.

The classical question of D.~H.~Lehmer 
(see~\cite[Problem~F12]{Guy}) is about the 
joint distribution of the parity of $n$ and 
its modular inverse $\inv{n}\in \cU_q$, defined 
by $n \inv{n} \equiv 1 \pmod q$.

W.~Zhang~\cite{Zha1,Zha2} has shown that the Weil bound
of Kloosterman sums, see~\cite[Corollary~11.12]{IwKow}, 
combined with some standard arguments, 
implies that if $q$ is odd then $n$ and 
its modular inverse $\inv{n}$ are of the same parity
$0.5 \varphi(q) + O\(q^{1/2 + o(1)}\)$ times for $n \in \cU_q$.

This result has been extended in generalised in 
various directions, including its multidimensional
analogues, 
see~\cite{ASZ,CoZa,LiuZha,LRS,Shp2,XuZha,YiZha,ZXY,ZhaZha,Zha1,Zha2,Zha3}
and references therein. 

In particular, it has been shown by Y.~Yi and W.~Zhang~\cite{YiZha},
that for any fixed integer $k\ne 0$, the smallest positive residue 
modulo $q$ of  $n^k$ and 
its modular inverse $\inv{n^k}$ are of the same parity
$0.5 \varphi(q) + O\(q^{3/4 + o(1)}\)$ times for $n \in \cU_q$.

Here we show that using the bound from~\cite{Shp1} of exponential sums
with sparse rational functions, one can get the same
(as in the case $k=1$) error  term
$O\(q^{1/2 + o(1)}\)$ for any fixed $k$ and in fact 
obtain an asymptotic formula in  a much more general 
case. Namely given an an integer $s\ge 2$ 
 $s$-dimensional integer vectors
$$
\vec{k} = (k_1, \ldots, k_s) , \qquad \vec{m} = (m_1, \ldots, m_s),
\qquad \vec{a} = (a_1, \ldots,a_s),
$$
where $m_1, \ldots, m_s \ge 1$, 
we denote by $N_q(\vec{m}, \vec{a}; \vec{k})$ the 
number of $n \in \cU_q$ such that  the smallest 
nonnegative residue of $n^{k_j}$ modulo $q$ is
 congruent to $a_j$ modulo $m_j$ 
for every $j =1, \ldots, s$.

In particular, the result of Y.~Yi and W.~Zhang~\cite{YiZha} can
be reformulated as the asymptotic formula,
\begin{equation}
\label{eq:YZ bound}
\begin{split}
N_q\((2,2), (0,0); (k,k)\) +   N_q\((2,2),   (1,1); (k,k)\)& \\
 = \frac{1}{2}  \varphi(q)~+~& O\(q^{3/4 + o(1)}\),
\end{split}
\end{equation}
which holds for any odd $q$.

Here we give the following generalisation and
improvement of~\eqref{eq:YZ bound}.

\begin{thm}
\label{thm:main}  
For any fixed integer $s\ge 2$ and 
a vector $\vec{k} \in \Z^s$ without zero components, 
uniformly over all vectors $\vec{a}, \vec{m} \in \Z^s$
with $m_1, \ldots, m_s \ge 1$ and an integer $q\ge 1$ with
$$
\gcd(m_1,q) =  \ldots = \gcd(m_s, q) = 1,
$$ 
we have 
$$
N_q(\vec{m}, \vec{a}; \vec{k})
= \frac{1}{m_1 \ldots m_s}  \varphi(q) + O\(q^{1-1/s + o(1)}\).
$$
\end{thm}

In particular, for $s=2$, $m_1 = m_2 = 2$ and 
$k_1 = k_2 = k$, Theorem~\ref{thm:main}  implies 
that the error term in the
asymptotic formula~\eqref{eq:YZ bound} is $O\(q^{1/2 + o(1)}\)$
for every odd $q\ge 1$.

Throughout the paper, the implied constants in the symbols `$O$',
and `$\ll$' may depend on the vector $\vec{k}$.
We recall that the notations $U = O(V)$ and $V \ll U$ are both 
equivalent to the assertion that the inequality $|U|\le cV$ holds for some
constant $c>0$.

\section{Exponential Sums}

For an integer $\ell$ we denote
$$
\el(z) = \exp(2 \pi i z/\ell)
$$
and recall that  for $u \in \Z$,
\begin{equation}
\label{eq:ident}
\frac{1}{\ell}
 \sum_{-(\ell-1)/2 \le \mu \le \ell/2}
 \el(\mu u)
=\left\{\begin{array}{rll}1 &\text{if}\ u
\equiv 0 \pmod \ell, \\ 0 &\text{otherwise,}\end{array}\right.
\end{equation} 
(which follows immediately from the formula for the
sum of a geometric progression).

We also recall  that for any integers  $U\ge 1$ and
$\mu$ with $0 < |\mu |\le  \ell/2$ we have
$$
\left|\sum_{u=0}^U \el(\mu u)\right| \ll 
\min\left\{U, \frac{\ell}{|\mu|}\right\},
$$
see~\cite[Bound~(8.6)]{IwKow}. 
In particular
\begin{equation}
\label{eq:Bound 1}
 \sum_{\substack{-(\ell-1)/2 \le \mu \le \ell/2\\ \mu \ne 0}}
\left|\sum_{u=0}^U \el(\mu u)\right| \ll \ell \log \ell
\end{equation} 
and
\begin{equation}
\label{eq:Bound 2}
 \sum_{-(\ell-1)/2 \le \mu \le \ell/2}
\left|\sum_{u=0}^U \el(\mu u)\right| \ll U + \ell \log \ell.
\end{equation}

Finally, as we have mentioned, our main tool
is the following slight generalisation 
of~\cite[Theorem~1]{Shp1} which gives an estimate 
of exponential sums with sparse rational functions.
We note that in~\cite[Theorem~1]{Shp1} only the case of $d=1$
has been considered, but the extension to
the case of arbitrary $d$ is straightforward.

\begin{lem}\label{lem:Sparse} For any fixed integer $s\ge 2$ and 
a vector $\vec{k} \in \Z^s$ without zero components  
and an integer $q\ge 1$, 
the bound
$$
\sum_{n \in \cU_q} 
\eq \( \sum_{1\le j\le s}\lambda_j n^{k_j} \)  \ll
d^{1/s}q^{1-1/s +o(1)}
$$
holds, uniformly over all integers $\lambda_1, \ldots, \lambda_s$
with 
$$
\gcd(\lambda_1, \ldots, \lambda_s) = d.
$$
\end{lem}

\begin{proof} Let
$$
q = \prod_{i=1}^\nu p_i^{\alpha_i}
$$ 
be the prime number factorization of $q$ and 
$q_i = q/p_i^{\alpha_i}$,  $i = 1, \ldots, \nu$.
We now define $t_i$ as the modular inverse of $q_i$ 
modulo $p_i^{\alpha_i}$, that is,
$$
t_iq_i \equiv \pmod {p_i^{\alpha_i}} \qquad \text{and}
\qquad 0 \le t_i < p_i^{\alpha_i},
$$
for $i = 1, \ldots, \nu$.
Using the multiplicative property of exponential sums
with rational functions, see~\cite[Equation~(12.21)]{IwKow}
or~\cite[Lemma~6]{Shp1}, we obtain 
\begin{equation}
\label{eq:CRT}
\sum_{n \in \cU_q} 
\eq \( \sum_{1\le j\le s}\lambda_j n^{k_j} \) 
=
\prod_{i=1}^\nu \sum_{n \in \cU_{p_i^{\alpha_i}}} 
\e_{p_i^{\alpha_i}} \( t_i\sum_{1\le j\le s}\lambda_j n^{k_j} \).
\end{equation}

Furthermore, by~\cite[Lemma~5]{Shp1} we have
\begin{equation}
\label{eq:Bound}
\sum_{n \in \cU_{p^{\alpha}}} 
\e_{p^\alpha} \(\sum_{1\le j\le s}\mu_j n^{k_j} \) \ll 
p^{\alpha(1 - 1/s + o(1))}
\end{equation}
for any prime power $p^\alpha$ with $p^{\alpha} \to \infty$
and integers $\mu_1, \ldots, \mu_s$ with 
$$
\gcd(\mu_1, \ldots, \mu_s,p) = 1.
$$
Combining~\eqref{eq:CRT} with~\eqref{eq:Bound} and
using that 
$$
\nu \ll \frac{\log q}{\log \log q}
$$
(which follows from the obvious inequality  
$\nu! \le q$ and the Stirling formula) 
we obtain the result.
\end{proof}

\subsection{Proof of Theorem~\ref{thm:main}}
\label{sec:discr}

Without loss of generality we may assume 
that $\vec{a}$ has a nonnegative components
satisfying $0 \le a_j < m_j$, $j =1, \ldots, s$. 

Let us define $U_j$ as the largest integers $U$
with $m_jU + a_j < q$, $j =1, \ldots, s$.  

Then $N_q(\vec{m}, \vec{a}; \vec{k})$ is equal to the 
number of solutions to the following system of congruences
\begin{equation}
\label{eq:Syst Prelim}
n^{k_j} \equiv m_j u_j + a_j \pmod q, 
\quad n \in \cU_q,\ 0 \le u_j \le U_j, \ 
j = 1, \ldots, s.
\end{equation}

Since $\gcd(m_1  \ldots  m_s, q) = 1$,  for every $j =1,\ldots, s$ 
we consider the moduluar inverse 
$r_j =\inv{m_j}$ of $m_j$ modulo $q$, 
and also define $b_j \in \cU_q$ by the congruence 
$b_j \equiv a_j r_j \pmod q$. 
Therefore, the system~\eqref{eq:Syst Prelim}
is equivalent to the 
the following system of congruences
\begin{equation}
\label{eq:Syst}
r_j n^{k_j} \equiv  u_j + b_j \pmod q, 
\quad n \in \cU_q,\ 0 \le u_j \le U_j, \ 
j = 1, \ldots, s.
\end{equation}

Using~\eqref{eq:ident} we write
\begin{eqnarray*}
N_q(\vec{m}, \vec{a}; \vec{k})
& = & 
\sum_{n \in \cU_q} \sum_{ 0 \le u_1 \le U_1}
\ldots \sum_{ 0 \le u_s \le U_s}\\
& & \qquad \frac{1}{q^s}
\sum_{-(q-1)/2 \le \lambda_1, \ldots, \lambda_s \le q/2}
\eq\(\sum_{1 \le j \le s} \lambda_j(r_j n^{k_j} - u_j - b_j)\).
\end{eqnarray*}
Changing the order of summation and then separating the
main term 
$$
\frac{\#\cU_q U_1\ldots U_s}{q^s} = 
\frac{\varphi(q) U_1\ldots U_s}{q^s}
$$
corresponding to $\lambda_1 = \ldots = \lambda_s = 0$, we
obtain
\begin{eqnarray*}
N_q(\vec{m}, \vec{a}; \vec{k}) & - & \frac{\varphi(q) U_1\ldots U_s}{q^s}\\
& = & \frac{1}{q^s}
\sum_{-(q-1)/2 \le \lambda_1, \ldots, \lambda_s \le q/2}
\hskip -45pt {\phantom{\sum}}^*  \hskip 25pt
\eq\(-\sum_{1 \le j \le s} \lambda_j b_j\)\\
& & \qquad \qquad \qquad\qquad 
\sum_{n \in \cU_q} \eq\(\sum_{1 \le j \le s} \lambda_j r_j n^{k_j} \)\\
& & \qquad \qquad \qquad \qquad \
\sum_{ 0 \le u_1 \le U_1}
\ldots \sum_{ 0 \le u_s \le U_s}\eq\(-\sum_{1 \le j \le s} \lambda_j u_j \).
\end{eqnarray*}
where $\Sigma^*$ means that the term corresponding to 
$\lambda_1 = \ldots = \lambda_s = 0$  is excluded from the summation. 
Therefore,
\begin{eqnarray*}
\lefteqn{
\left|N_q(\vec{m}, \vec{a}; \vec{k}) -
 \frac{\varphi(q) U_1\ldots U_s}{q^s}\right|}\\
&   & \qquad \le \frac{1}{q^s}
\sum_{-(q-1)/2 \le \lambda_1, \ldots, \lambda_s \le q/2}
\hskip -45pt {\phantom{\sum}}^*  \hskip 25pt  
\left|\sum_{n \in \cU_q}
 \eq\(\sum_{1 \le j \le s} \lambda_j r_j n^{k_j} \)\right|\\
& & \qquad \qquad \qquad \qquad \qquad \qquad \qquad 
\prod_{1 \le j \le s} 
\left|\sum_{ 0 \le u_j \le U_j} \eq\(\lambda_j u_j\)\right|.
\end{eqnarray*}

Now, for every divisor $d \mid q$ we collect together the terms 
with the same value of 
$\gcd(\lambda_1, \ldots, \lambda_s) = d$ and 
then apply Lemma~\ref{lem:Sparse}, obtaining the 
estimate
\begin{equation}
\label{eq:N and Sigmad}
\left|N_q(\vec{m}, \vec{a}; \vec{k}) -
 \frac{\varphi(q) U_1\ldots U_s}{q^s}\right|
\le  q^{-s + 1 - 1/s + o(1)}
\sum_{\substack{d\mid q\\ q < q}} d^{1/s} \Sigma_d,
\end{equation}
where
$$
\Sigma_d  =
\sum_{\substack{-(q-1)/2 \le \lambda_1, \ldots, \lambda_s \le q/2
\\ \gcd(\lambda_1, \ldots, \lambda_s) = d}} \ \prod_{1 \le j \le s} 
\left|\sum_{ 0 \le u_j \le U_j} \eq\(\lambda_j u_j\)\right|.
$$
Writing $\lambda_j = d \mu_j$, $j=1, \ldots, s$, and $q = d q_d$, 
we derive
$$
\Sigma_d  =
\sum_{\substack{-(q_d-1)/2 \le \mu_1, \ldots, \mu_s \le q_d/2
\\ \gcd(\mu_1, \ldots, \mu_s) = 1}} \ \prod_{1 \le j \le s} 
\left|\sum_{ 0 \le u_j \le U_j} \eqd\(\mu_j u_j\)\right|.
$$
Furthermore, we have
\begin{equation}
\label{eq:Sigmadj}
\Sigma_d  \le  \sum_{1 \le j \le s} \sigma_{d,j}
\end{equation}
where 
$$
\sigma_{d,j} = \sum_{\substack{-(q_d-1)/2 \le \mu_1, \ldots, \mu_s \le q_d/2
\\ \mu_j  \ne 0}} \ \prod_{1 \le j \le s} 
\left|\sum_{ 0 \le u_j \le U_j} \eqd\(\mu_j u_j\)\right|,
\qquad j =1, \ldots, s.
$$
We have,
\begin{eqnarray*}
\sigma_{d,j} & = & \sum_{\substack{-(q_d-1)/2 \le \mu_j \le q_d/2
\\ \mu_j  \ne 0}} \left|\sum_{ 0 \le u_j \le U_j} \eqd\(-\mu_j u_j\)\right|\\
& & \qquad\qquad \prod_{\substack{1 \le h \le s\\  h \ne j}} \ 
\sum_{-(q_d-1)/2 \le \mu_h \le q_d/2}
\left|\sum_{ 0 \le u_h \le U_h} \eqd\(\mu_h u_h\)\right|. 
\end{eqnarray*}
Applying~\eqref{eq:Bound 1} to the sum over $\mu_j$  and~\eqref{eq:Bound 2} 
to the other sums (and using the trivial
estimate $U_h \le q$), we obtain
$$
\sigma_{d,j} \le (q_d \log q_d) \prod_{\substack{1 \le h \le s\\  h \ne j}}
(U_h + q_d \log q_q)  \le q_d q^{s-1} (\log q_d)^s = 
d^{-1} q^{s + o(1)}
$$
for every $j=1, \ldots, s$. Substituting this in~\eqref{eq:Sigmadj},
and then recalling~\eqref{eq:N and Sigmad}, we obtain
$$
\left|N_q(\vec{m}, \vec{a}; \vec{k}) -
 \frac{\varphi(q) U_1\ldots U_s}{q^s}\right|
\le  q^{1 - 1/s + o(1)}
\sum_{\substack{d\mid q\\ q < q}} d^{-1 + 1/s}
\le    q^{1 - 1/s + o(1)} 
\sum_{d\mid q} 1.
$$
By the well-known estimate on the 
divisor function,
$$
\sum_{d\mid q} 1 = q^{o(1)}
$$ 
see~\cite[Theorem~317]{HardyWright}, we obtain
$$
\left|N_q(\vec{m}, \vec{a}; \vec{k}) -
 \frac{\varphi(q) U_1\ldots U_s}{q^s}\right|
\le  q^{1 - 1/s + o(1)}.
$$
It remains to notice that $U_j = q/m_j + O(1)$,  $j=1, \ldots, s$.
Therefore 
$$
U_1\ldots U_s = \frac{q^s}{m_1\ldots m_s} + O(q^{s-1}), 
$$
which concludes the proof.


\begin{thebibliography}{100}

\bibitem{ASZ} E.  Alkan, F.~Stan and A. Zaharescu, 
`Lehmer $k$-tuples', 
{\it Proc. Amer. Math. Soc.\/}, {\bf 134} (2006), 2807--2815. 
 
 \bibitem{CoZa} C. Cobeli and A. Zaharescu, `Generalization of a 
problem of Lehmer', 
{\it Manuscr. Math.\/}, {\bf 104} (2001),  301--307.
 
\bibitem{Guy} R. K. Guy, {\it Unsolved problems in number theory\/},
Springer-Verlag, New York, 1994.

\bibitem{HardyWright} G. H. Hardy and E. M. Wright, {\it An introduction to
the theory of numbers\/}, Oxford Univ. Press, Oxford, 1979.

\bibitem{IwKow} H. Iwaniec and E. Kowalski,  {\it Analytic number
theory\/},  Amer. Math. Soc., Providence, RI, 2004.

\bibitem{LiuZha} H. N.  Liu and W. Zhang,`On a
problem of D.~H.~Lehmer', {\it Acta Math. Sinica\/},
{\bf 22} (2006), 61--68.

\bibitem{LRS} S. R. Louboutin, J. Rivat and  A. S{\'a}rk{\"o}zy',
`On a problem of D. H. Lehmer',
{\it Proc. Amer. Math. Soc.\/},  {\bf 135} (2007),
969--975.

\bibitem{Shp1} I. E. Shparlinski, 
`On exponential sums with sparse
polynomials and rational functions',  
{\it J. Number  Theory\/}, {\bf 60} (1996), 233--244.


\bibitem{Shp2} I. E. Shparlinski, `On a generalisation of a Lehmer
problem', {\it Preprint\/}, 2006 (available from {\tt
http://arxiv.org/abs/math/0607414}). 

 
 \bibitem {XuZha} Z. Xu and W. Zhang,
`On a problem of D.~H.~Lehmer over short intervals',
{\it J. Math. Anal. Appl.\/}, {\bf 320} (2006), 756--770.
 
\bibitem {YiZha}  Y. Yi and W. Zhang,
`On the generalization of a problem of D. H. Lehmer',
{\it Kyushu J. Math.\/}, {\bf  56}  (2002), 235--241.

\bibitem {ZXY} W. Zhang, Z. Xu  and Y. Yi, `A problem of D. H. Lehmer
and its mean square value formula',
{\it J. Number Theory\/}, {\bf 103} (2003),  197--213.

\bibitem {ZhaZha}  T. Zhang and W. Zhang, 
`A generalization on the difference between an integer
and its inverse modulo  $q$, II',  
{\it Proc. Japan Acad. Sci., Ser.A\/}, {\bf 81} (2005),  7--11. 

\bibitem {Zha1}  W. Zhang, 
`On a problem of D.~H.~Lehmer and its generalization',  
{\it Compos. Math.\/}, {\bf 86} (1993),  307--316. 

\bibitem {Zha2}  W. Zhang, 
`On a problem of D.~H.~Lehmer and its generalization, II',  
{\it Compos. Math.\/}, {\bf 91} (1994),   47--56. 

%\bibitem {Zha3}  W. Zhang, 
%`On the difference between a D.~H.~Lehmer number and its 
%inverse modulo  $q$',  
%{\it Acta Arith.\/}, {\bf 68} (1994),   255--263.


\bibitem {Zha3}
W. Zhang, `On a problem of D.~H.~Lehmer and Kloosterman sums',
{\it Monatsh. Math.\/}, {\bf 139} (2003), 247--257.


\end{thebibliography}
\end{document}